\documentclass[11pt]{article}
\usepackage{latexsym,amsfonts,amssymb}
\input amssym.def

\begin{document}
\baselineskip 22pt
\title{\bf The classification of $\bf Z\!-\!$graded modules of the intermediate series
over the $q$-analog Virasoro-like algebra$^{\ast}$
\thanks{ $\ast$ Supported by the National Science Foundation of China
(No. 10671160) and the China Postdoctoral Science Foundation (No.
20060390693). To appear in Algebra colloquium.}} \vskip1cm
\author{\small Yina Wu$^1$ and Weiqiang Lin$^{1,2}$
\thanks{e-mail: linwq83@yahoo.com.cn}\\
 \small 1. Department of Mathematics, Zhangzhou Teacher's College,\\
\small Zhangzhou  363000, Fujian, China.\\
\small 2. Department of Mathematics, University of Science and Technology of China,\\
\small Hefei\ 230026, Anhui, China.
 }
\date{}
\maketitle{{\noindent\bf Abstract}\\
 \indent In this paper, we complete the classification of the {\bf Z}-graded modules of the intermediate series over the
 $q$-analog Virasoro-like algebra $L$. We first
construct four classes of irreducible {\bf Z}-graded $L$-modules of
the intermediate series.
 Then we prove that any {\bf Z}-graded $L$-modules of the intermediate series must
 be
 the direct sum of some trivial $L$-modules or one of the modules constructed by us.

{\noindent\small{\it Keywords:}  modules of the intermediate series,
{\bf Z}-graded $L$-module, the $q$-analog Virasoro-like
algebra.}

}

\section{Introduction}
\indent $\ \ \ \ \ $The classification of irreducible graded modules
with finite dimensional homogeneous subspaces over a graded Lie
algebra is one of the main subject in the study of Lie theory.
Meanwhile, the irreducible graded modules with finite dimensional
homogeneous subspaces for some infinite dimensional Lie algebras,
such as the Heisenberg Lie algebra and the Virasoro algebra, have
important applications in the study of the vertex operator algebras
and theoretical physics. In this paper, we study the classification
of {\bf Z}-graded modules of the intermediate series over the
$q$-analog Virasoro-like algebra. The $q$-analog Virasoro-like
algebra is introduced by Kirkman etc in [1]. It can be realized as
the universal central extension of the inner derivation Lie algebra
of the quantum torus ${\bf C}_{q}[x_{1}^{\pm1},x_{2}^{\pm1}]$ (see
[2] or [3]), where $q$ is generic. Quantum torus is one of the main
objects in noncommutative geometry, and plays an important role in
the classification of extended affine Lie algebras ([2]). Meanwhile,
the $q$-analog Virasoro-like algebra can be regarded as a $q$
deformation of the Virasoro-like algebra introduced and studied by
Arnold, Wit, etc when they tried to generalize the Virasoro algebra
([4], [5] and [1]). There are some papers devoted to the study of
the structure and representations of the $q$-analog Virasoro
algebra. Jiang and Meng studied its derivation Lie algebra and the
automorphism group of its derivation Lie algebra ([3]). Chen, Lin,
etc. studied the structure of its automorphism group ([6]). Zhao and
Rao constructed a class of highest weight irreducible {\bf Z}-graded
modules over the $q$-analog Virasoro-like algebra, and gave a
sufficient and necessary condition for such a module with finite
dimensional homogeneous subspaces ([7]). Gao constructed a class of
principal vertex representations for the extended affine Lie
algebras coordinatized by certain quantum tori by using the
representation of the $q$-analog Virasoro-like algebra([8]). The
classification of the irreducible graded modules with finite
dimensional homogeneous subspaces and nontrivial centers over the
$q$-analog Virasoro-like algebra has been completed in [9]. Thus, we
only consider the classification of the {\bf Z}-graded modules of
the intermediate series over the centerless $q$-analog Virasoro-like
algebra $L$ in the present paper. In section 2, we first recall some
notations about the centerless $q$-analog Virasoro-like algebra $L$
and its $\bf Z$-graded modules of the intermediate series. Then we
construct four classes of {\bf Z}-graded $L$-modules of the
intermediate series, and show that they are irreducible. In section
3, we complete the classification of the {\bf Z}-graded $L$-modules
of the intermediate series.
\par

\section{The $q$-analog Virasoro-like algebra and its {\bf Z}-graded
modules of the intermediate series $\hspace*{\fill}$} \indent $\ \ \
\ \ $ Throughout this paper we use $\bf {Z,\ Z^{*},\ C,\ C^{*}}$ and
{\bf N} to denote the sets of integers, nonzero integers, complex
number, nonzero complex number and positive integers respectively.
All spaces are over {\bf C}.\par In this paper, we require $q\in \bf
C$ to be a fixed nonzero non-root of unity. Let $L$ be a vector
space spanned by $\left\{t_{1}^{h}t_{2}^{j}| (h,j)\in {\bf
Z}^{2}\setminus \{(0,0)\}\right\}$. We denote it by $L=\left\langle
t_{1}^{h}t_{2}^{j}\ |\ (h,j)\in {\bf Z}^{2}\setminus
\{(0,0)\}\right\rangle$. Define a commutator in $L$ as follows:
$$[t_{1}^{h}t_{2}^{j},t_{1}^{m}t_{2}^{n}]=\left(q^{jm}-q^{hn}\right)t_{1}^{h+m}t_{2}^{j+n},\eqno(2.1)$$
then $L$ is called the centerless $q$-analog Virasoro-like algebra.
And it is easy to check that $L=\bigoplus\limits_{u\in {\bf
Z}}L_{u}$, where $L_{u}=\left\langle t_{1}^{u}t_{2}^{j}|j\in {\bf
Z}\right\rangle$, is a {\bf Z}-graded Lie algebra.
\par Next, we recall the definition of the {\bf Z}-graded $L$-modules of the intermediate series.
If a vector
space $V=\bigoplus\limits_{n\in {\bf Z}}V_{n}$ satisfies:\\
(1) $V$ is a $L$-module,\\
(2) Regarding to
$L=\bigoplus\limits_{u\in {\bf Z}}L_{u}$, there are
$L_{u}.V_{n}\subseteq V_{n+u}\ \mbox{for\ any}\ n,\ u\in {\bf Z},$\\
(3) $dim V_{n}\leq 1,\ \forall\ n\in {\bf Z},$\\
then $V=\bigoplus\limits_{n\in {\bf Z}}V_{n}$ is called a {\bf
Z}-graded $L$-module of the intermediate series.\par \par In this
paper, we study this type of {\bf Z}-graded $L$-module and its
classification. First we construct four classes of {\bf Z}-graded
$L$-modules.\par
 {\bf Proposition 2.1}$:$  Set V=$\langle v_{j}|j\in {\bf
 Z}\rangle$. For any $a\in {\bf C}^{\ast}$, define the action of the elements in $L$ on $V$ by
 linearly extending the following maps
 respectively:\\
 $\langle1\rangle.\
(t_{1}^{m}t_{2}^{n}).v_{k}=(aq^{k})^{n}v_{k+m},\ for\ all\ (m,n)\in
{\bf Z}^{2}\setminus\{(0,0)\};\\
\langle2\rangle.\
(t_{1}^{m}t_{2}^{n}).v_{k}=(-1)^{m}(aq^{k})^{n}v_{k+m},\ for\ all\
(m,n)\in {\bf Z}^{2}\setminus\{(0,0)\};\\
\langle3\rangle.\
(t_{1}^{m}t_{2}^{n}).v_{k}=(-1)^{m+n+1}(aq^{-k-m})^{n}v_{k+m},\
for\ all\ (m,n)\in {\bf Z}^{2}\setminus\{(0,0)\};\\
\langle4\rangle.\
(t_{1}^{m}t_{2}^{n}).v_{k}=(-1)^{n+1}(aq^{-k-m})^{n}v_{k+m},\ for\
all\ (m,n)\in {\bf Z}^{2}\setminus\{(0,0)\}.$\\
Then $V$ becomes a $L$-module with the actions defined in
$\langle1\rangle$, $\langle2\rangle$, $\langle3\rangle$ or
$\langle4\rangle$ respectively.
\par  {\bf Proof.}
We take $\langle1\rangle\ and\  \langle3\rangle$ for examples to
show that $V$ is a $L$-module with respect to the action defined in
$\langle1\rangle \mbox{or} \langle3\rangle$ respectively.
\par For
the operators defined in $\langle1\rangle$, we have
$$
[t_{1}^{m}t_{2}^{n},t_{1}^{j}t_{2}^{h}].v_{k}=(q^{jn}-q^{hm})(t_{1}^{m+j}t_{2}^{n+h}).v_{k}
=(q^{jn}-q^{hm})(aq^{k})^{n+h}v_{k+m+j}.
$$
Meanwhile,
$$
(t_{1}^{m}t_{2}^{n}).(t_{1}^{j}t_{2}^{h}).v_{k}-(t_{1}^{j}t_{2}^{h}).(t_{1}^{m}t_{2}^{n}).v_{k}
=(aq^{k})^{h}(t_{1}^{m}t_{2}^{n}).v_{k+j}-(aq^{k})^{n}(t_{1}^{j}t_{2}^{h}).v_{k+m}
$$
$$
=(aq^{k})^{h}(aq^{k+j})^{n}v_{k+j+m}-(aq^{k})^{n}(aq^{k+m})^{h}v_{k+j+m}$$
$$=(q^{jn}-q^{hm})(aq^{k})^{n+h}v_{k+m+j}.$$
Hence,\ $$[t_{1}^{m}t_{2}^{n},t_{1}^{j}t_{2}^{h}].v_{k}=
(t_{1}^{m}t_{2}^{n}).(t_{1}^{j}t_{2}^{h}).v_{k}-(t_{1}^{j}t_{2}^{h}).(t_{1}^{m}t_{2}^{n}).v_{k}\
.$$
 Therefore,\ $V$\ is\ a\ $L$-module.\par For the operators defined
 in
$\langle3\rangle$, we have
$$[t_{1}^{m}t_{2}^{n},t_{1}^{j}t_{2}^{h}].v_{k}=(q^{jn}-q^{hm})(t_{1}^{m+j}t_{2}^{n+h}).v_{k}
=(q^{jn}-q^{hm})(-1)^{m+j+n+h+1}(aq^{-k-m-j})^{n+h}v_{k+m+j},$$
and
$$(t_{1}^{m}t_{2}^{n}).(t_{1}^{j}t_{2}^{h}).v_{k}-(t_{1}^{j}t_{2}^{h}).(t_{1}^{m}t_{2}^{n}).v_{k}$$
$$=(-1)^{j+h+1}(aq^{-k-j})^{h}(t_{1}^{m}t_{2}^{n}).v_{k+j}-
(-1)^{m+n+1}(aq^{-k-m})^{n}(t_{1}^{j}t_{2}^{h}).v_{k+m}$$
$$=(q^{jn}-q^{hm})(-1)^{m+j+n+h+1}(aq^{-k-m-j})^{n+h}v_{k+m+j}
.$$ Thus
$$[t_{1}^{m}t_{2}^{n},t_{1}^{j}t_{2}^{h}].v_{k}=
(t_{1}^{m}t_{2}^{n}).(t_{1}^{j}t_{2}^{h}).v_{k}-(t_{1}^{j}t_{2}^{h}).(t_{1}^{m}t_{2}^{n}).v_{k}
.$$ Hence, $V$ is a $L$-module.\par Similarly, one can check that
$V$ is a $L$-module with the operators defined in $\langle2\rangle$
or $\langle4\rangle$ respectively. $\hspace*{\fill} \blacksquare$

   {\bf Remark 2.2}$:$  One can easily see that the modules defined above are irreducible {\bf Z}-graded
   $L$-modules of the intermediate series with respect to the linear space decomposition $V=\bigoplus\limits_{n\in {\bf
   Z}}V_{n}$, where $V_{n}={\bf C}v_{n}$. We will use $V(a,I),\ V(a,II),\ V(a,III)$ and $V(a,IV)$ to denote the
   corresponding {\bf Z} graded
$L$-modules defined by $\langle1\rangle,\ \langle2\rangle,\
\langle3\rangle,\ \mbox{and}\ \langle4\rangle$ respectively.

\section{The classification of Z-graded $L$-modules of the intermediate series}
\indent In this section, we discuss the classification of the {\bf
Z}-graded $L$-modules of the intermediate series. We will first
prove two Lemmas.
\par {\bf Lemma 3.1}$:$ If $V=\bigoplus\limits_{n\in {\bf
Z}}V_{n}$ is a {\bf Z}-graded $L$-modules of the intermediate series
and the action of $t_{1}^{1}.t_{1}^{-1}$ is degenerate, then $V$ can
be decomposed into the direct sum of some trivial $L$-submodules.
\par { \bf Proof}.
Without loss of generality, we can assume $V\neq0$. Since the action
of $t_{1}^{1}.t_{1}^{-1}$ is degenerate, there exists a nonzero
vector $v_{j}\in V_j$ such that $t_{1}^{1}.t_{1}^{-1}.v_{j}=0$. Thus
$t_{1}^{1}.t_{1}^{-1}.v_{j}=t_{1}^{-1}.t_{1}^{1}.v_{j}=0$ since
$[t_1^1,t_1^{-1}]=0$ by the definition of the Lie algebra $L$.
Therefore, we obtain that
$$\left\{\begin{array}{l} t_{1}^{1}.V_{j}=0;\\t_{1}^{-1}.V_{j}=0,\end{array}\right.
\left\{\begin{array}{l}
t_{1}^{-1}.V_{j}=0;\\t_{1}^{-1}.V_{j+1}=0,\end{array}\right.
\left\{\begin{array}{l}
t_{1}^{1}.V_{j-1}=0;\\t_{1}^{1}.V_{j}=0,\end{array}\right.
\mbox{or}\
 \left\{\begin{array}{l}
t_{1}^{1}.V_{j-1}=0;\\t_{1}^{-1}.V_{j+1}=0,\end{array}\right.\eqno(3.1)$$
since dim $V_{k}\leq1$ for all $k\in {\bf Z}$. Considering that
$t_{1}^{1}.t_{1}^{-1}.v=t_{1}^{-1}.t_{1}^{1}.v$ for any $v\in V$ and
dim $V_{n}\leq1,\forall\ n\in {\bf Z}$, we deduce that, for any
$n\in{\bf Z}$,
$$
(1)\left\{\begin{array}{l}
t_{1}^{1}.V_{n}=0;\\t_{1}^{-1}.V_{n}=0.\end{array}\right.
(2)\left\{\begin{array}{l}
t_{1}^{-1}.V_{n}=0;\\t_{1}^{-1}.V_{n+1}=0.\end{array}\right.
(3)\left\{\begin{array}{l}
t_{1}^{1}.V_{n-1}=0;\\t_{1}^{1}.V_{n}=0.\end{array}\right.
\mbox{or}\
(4)\left\{\begin{array}{l}
t_{1}^{1}.V_{n-1}=0;\\t_{1}^{-1}.V_{n+1}=0.\end{array}\right.\eqno(3.2)$$
\par Now we prove that $V_n$ is a trivial $L$-module. We will first prove that $t_{2}^{k}$ acts trivially on
$V_{n}$ for any $k\in {\bf Z}^{*},\ n\in{\bf Z}$. We only give the
proof of this claim for case (1) here. The proofs of this claim for
the other three cases are similar.
\par  Suppose that $t_{1}^{1}.V_{n}=0$ and
$t_{1}^{-1}.V_{n}=0$. By the definition of {\bf Z}-graded
$L$-module, we have that
$$ t_{2}^{k}.V_{n}\subseteq V_{n}\ \mbox{for\ all}\ k\in
{\bf Z}^{*}.\eqno(3.3)$$ Thus we deduce that
$$[[t_{1}^{1},t_{2}^{k}],t_{1}^{-1}].V_{n}=(t_{1}^{1}.t_{2}^{k}.t_{1}^{-1}.V_{n}-
t_{2}^{k}.t_{1}^{1}.t_{1}^{-1}.V_{n})-t_{1}^{-1}.(t_{1}^{1}.t_{2}^{k}.V_{n}-
t_{2}^{k}.t_{1}^{1}.V_{n})=0.$$
Therefore $t_{2}^{k}$ acts trivially
on $V_{n}$ since
$$[[t_{1}^{1},t_{2}^{k}],t_{1}^{-1}].V_{n}=(1-q^{k})(q^{-k}-1)t_{2}^{k}.V_{n}.$$
\par Next, we show that $t_1^j$ acts trivially on $V_n$ for any
$j\in {\bf Z}^{\ast}$.
 On one hand, we have that
$$[[t_{1}^{j},t_{2}^{1}],t_{2}^{-1}].V_{n}=(1-q^{j})(1-q^{-j})t_{1}^{j}.V_{n}.$$
On the other hand, for any $j\in {\bf Z}^{*}$, we have that
$$[[t_{1}^{j},t_{2}^{1}],t_{2}^{-1}].V_{n}=(t_{1}^{j}.t_{2}^{1}.t_{2}^{-1}.V_{n}-
t_{2}^{1}.t_{1}^{j}.t_{2}^{-1}.V_{n})-t_{2}^{-1}.(t_{1}^{j}.t_{2}^{1}.V_{n}-
t_{2}^{1}.t_{1}^{j}.V_{n})=0,$$ since $t_{2}^{k}$ acts trivially on
$V_{n},\ \forall \ k\in{\bf Z}^*, \ n\in{\bf Z}$.  Hence,
$t_{1}^{j}.V_{n}=0$ for any $j\in {\bf Z}^{*}$.
\par Finally, for
any $k,j\in{\bf Z}^*$, we have that
$$
(1-q^{jk})t_1^jt_2^k.V_n=[t_1^j,t_2^k].V_n=t_1^j.t_2^k.V_n-t_2^k.t_1^j.V_n=0.
$$
Therefore, $V_{n}$ is a trivial L-submodules for any $n\in {\bf Z}$.
Thus $V$ can be decomposed into the direct sum of some trivial
$L$-submodules. $\hspace*{\fill} \blacksquare$

 \par {\bf Lemma 3.2}$:$
 If $V=\bigoplus\limits_{n\in {\bf Z}}V_{n}$ is a {\bf Z}-graded
$L$-modules of the intermediate series and the action of
$t_{1}^{1}.t_{1}^{-1}$ is nondegenerate, then V must be isomorphic
to $V(a,I),\ V(a,II),\ V(a,III)$ or $V(a, IV)$ for some $a\in{\bf
C}^*$. (ref. {\bf Remark 2.2})
\par { \bf Proof}.
Since the action of $t_{1}^{1}.t_{1}^{-1}$ is nondegenerate, the
action of $t_{1}^{\pm1}$ is nondegenerate. Together with
$L_{u}.V_{n}\subseteq V_{n+u}$ and dim $V_{n}\leq 1$ for any $u\in
{\bf Z}^{*}\ and\ n\in {\bf Z}$, there must be
$t_{1}^{\pm1}.V_{n}=V_{n\pm1}$. Thus dim $V_{n}$=1 for any $n\in
{\bf Z}.$\par

We first show that there exists a base $\{v_{j}\in V_{j}| j\in {\bf
Z}\}$ of V  such that $t_{1}^{\pm 1}. v_{j}=\lambda v_{j\pm1}$ for
any $j\in {\bf Z}$. Suppose $\omega_{0}\in V_{0}$ with
$\omega_{0}\neq 0$ and set $\omega_{n}=t_{1}^{1}.\omega_{n-1}\in
V_{n}$. Since the action of $t_{1}^{\pm 1}$ is nondegenerate and dim
$V_{n}=1$ for all $n\in {\bf Z}$, we have that $\{\omega_{n}\in
V_{n}|n\in {\bf Z}\}$ forms a base of $V$ and
$t_{1}.\omega_{n}=\omega_{n+1}$. Denote $t_{1}^{-1}.
\omega_{k}=\phi(k)\omega_{k-1}$. By (2.1), we have
$[t_{1}^{1},t_{1}^{-1}]=0$. Thus $t_{1}.(t_{1}^{-1}.
\omega_{k})=t_{1}^{-1}.(t_{1}.\omega_{k})$ for all $k\in {\bf Z}$,
which implies $\phi(k)=\phi(k+1)$ for any  $k\in {\bf Z}$. Thus
there exists $p\in{\bf C}^*$ such that $t_{1}^{-1}.\omega_{k}=p
\omega_{k-1}$ for any $k\in {\bf Z}$. Set $\lambda=\sqrt{p}\neq 0$
and $v_{k}=\frac{\omega_{k}}{\lambda^{k}}$, then $t_{1}^{\pm1}.
v_{k}=\lambda v_{k\pm1}$ for any $k\in {\bf Z}$.

\par Set $(t_{1}^{h}t_{2}^{j}). v_{k}=f(h,j,k)v_{k+h}$.
 Then $f(\pm1,0,k)=\lambda$ for any $k\in {\bf
Z}$ by the choice of $v_k$. Now, we prove that f(m,0,k)=f(m,0,0) for
any $k\in {\bf Z}$. Since $[t_{1}^{m},t_{1}^{1}]=0$, we  have that
$t_{1}^{m}.(t_{1}^{1}. v_{k})=t_{1}^{1}.(t_{1}^{m}. v_{k})$ which
implies
$$f(m,0,k)=f(m,0,k+1)\ \mbox{for\ any}\ m\in {\bf Z}^{*},\ k\in {\bf
Z}.$$ Therefore, $$f(m,0,k)=f(m,0,0)\ \mbox{for\ any}\ m\in {\bf
Z}^{*},\ k\in {\bf Z}. \eqno(3.4)$$
\par Next we prove the following claim.

{\bf Claim:} $f(0,1,k)=ab^{k},\
f(0,-1,k)=\frac{(1-q)(1-q^{-1})}{a(2-b-b^{-1})}b^{-k}$
 for some $a,\ b\in {\bf
C}^{*}\setminus\{1\}$.
\par Since $[[t_{1}^{1},t_{2}^{\pm1}],t_{1}^{-1}].
v_{k}=(1-q)(q^{-1}-1)t_{2}^{\pm1}. v_{k}$, we obtain that
$$
(2\lambda^{2}+(1-q)(q^{-1}-1))f(0,\pm1,k)=\lambda^{2}(f(0,\pm1,k-1)+f(0,\pm1,k+1)),\
\forall\ k\in {\bf Z}.\eqno (3.5)
$$
We have that the characteristic equation of (3.5) is as follow
$$\lambda^{2}x^{2}-(2\lambda^{2}+(1-q)(q^{-1}-1))x+\lambda^{2}=0.\eqno(3.6)$$
\par Now we divide our proof of the claim into two cases according to whether the
equation (3.6) has different roots or not.
\par Case I. The equation (3.6) has not different roots. Then
$$\triangle=
4\lambda^{2}(1-q)(q^{-1}-1)+(1-q)^{2}(q^{-1}-1)^{2}=0,
$$
which implies $(1-q)(q^{-1}-1)=-4\lambda^{2}$. Applying this result
to (3.6), we have $x^{2}+2x+1=0$, since $\lambda\neq0$. Thus the
root of equation (3.6) is $-1$. By a result in [10], we obtain that
there exist $\lambda_1(\pm 1),\lambda_2(\pm 1)\in {\bf C}$ such that
$$ f(0,\pm1,k)=(-1)^{k}\left(\lambda_{1}(\pm1))+k\lambda_{2}(\pm1)\right),\
\forall \ k\in {\bf Z}.\eqno (3.7)$$ Since
$[[t_{1}^{1},t_{2}^{1}],t_{2}^{-1}].v_{k}=(1-q)(1-q^{-1})t_{1}^{1}.v_{k}$,
 we have
$$(1-q)(1-q^{-1})=(f(0,1,k)-f(0,1,k+1))(f(0,-1,k)-f(0,-1,k+1)).\eqno(3.8)$$
Applying (3.7) to (3.8), we have that
$$(1-q)(1-q^{-1})=4\lambda_{1}(1)\lambda_{1}(-1)+
(2k+1)\left(2\lambda_{2}(1)\lambda_{1}(-1)+2\lambda_{1}(1)\lambda_{2}(-1)+
(2k+1)\lambda_{2}(1)\lambda_{2}(-1)\right)$$ for all $k\in {\bf Z}$,
which implies that
$$\lambda_{2}(1)\lambda_{2}(-1)=0,\;\;
\lambda_{2}(1)\lambda_{1}(-1)+\lambda_{1}(1)\lambda_{2}(-1)=0,
$$ and
$$
4\lambda_{1}(1)\lambda_{1}(-1)=(1-q)(1-q^{-1}). $$
Therefore
$\lambda_{2}(1)\lambda_{1}(-1)=\lambda_{1}(1)\lambda_{2}(-1)=0.$
Thus we obtain that
$$\lambda_{2}(1)=\lambda_{2}(-1)=0\ \mbox{but}\
\lambda_{1}(1)\lambda_{1}(-1)\neq0.$$
\par Hence $f(0,1,k)=\lambda_{1}(1)(-1)^{k}\ and\
\ f(0,-1,k)=\frac{(1-q)(1-q^{-1})}{4\lambda_{1}(1)}(-1)^{k}$. Thus
the claim holds in this case with $a=\lambda(1)$ and $b=-1$.

\par Case II. The equation (3.6) has different roots. Then
$$\triangle:=4\lambda^{2}(1-q)(q^{-1}-1)+(1-q)^{2}(q^{-1}-1)^{2}\neq0.$$
Since $\lambda\neq0$, the root of equation (3.6) can not be zero.
Thus we can assume the roots of (3.6) are $x$ and $x^{-1}$
respectively. By a result in [10], we have that there exist
$\lambda_1(\pm 1),\lambda_2(\pm1)\in{\bf C}$ such that
$$f(0,\pm1,k)=\lambda_{1}(\pm1)x^{k}+\lambda_{2}(\pm1)x^{-k},\
\forall\ k\in {\bf Z}.\eqno(3.9)$$ Substituting (3.9) into (3.8), we
obtain that
$$(1-q)(1-q^{-1})=\lambda_{1}(-1)\lambda_{1}(1)(1-2x+x^{2})x^{2k}+
\lambda_{2}(-1)\lambda_{2}(1)(1-2x^{-1}+x^{-2})x^{-2k}$$
$$-\left(\lambda_{1}(-1)\lambda_{2}(1)
+\lambda_{2}(-1)\lambda_{1}(1)\right)(x+x^{-1}-2),\eqno(3.10)$$ for
all $k\in {\bf Z}$.
\par Now we prove that
$$\lambda_{1}(-1)\lambda_{1}(1)=\lambda_{2}(-1)\lambda_{2}(1)=0.$$
 We
will divide our proof into two cases according to whether $|x|=1$ or
not.

\par Subcase 1. $|x|=1$.
\par Since the equation (3.6) has different roots, we deduce that $x\neq \pm1$.
If $x^{2}\neq -1$, as the equation (3.10) holds for all $k\in {\bf
Z}$, then it is easy to obtain
$\lambda_{1}(-1)\lambda_{1}(1)=\lambda_{2}(-1)\lambda_{2}(1)=0$ by
the geometry significance. If $x^{2}=-1$, then $x=\pm i$. Without
loss of generality, we can assume $x=i$. Substituting it into (3.9),
we have
$$f(0,\pm1,k)=\lambda_{1}(\pm1)i^{k}+\lambda_{2}(\pm1)(-i)^{k}.\eqno(3.11)$$
Substituting x=i into (3.10), we have
$$(1-q)(1-q^{-1})=(-1)^{k}(2i)\left(-\lambda_{1}(-1)\lambda_{1}(1)
+\lambda_{2}(-1)\lambda_{2}(1)\right)-2
\left(\lambda_{1}(-1)\lambda_{2}(1)+\lambda_{2}(-1)\lambda_{1}(1)\right),\eqno(3.12)$$
where $k\in {\bf Z}$. Notice that (3.12) holds for any $k\in {\bf
Z}$, we deduce that
$$\lambda_{1}(-1)\lambda_{1}(1)=\lambda_{2}(-1)\lambda_{2}(1). $$
If
$\lambda_{1}(-1)\lambda_{1}(1)=\lambda_{2}(-1)\lambda_{2}(1)\not=0$,
we can assume
$\frac{\lambda_{1}(1)}{\lambda_{2}(1)}=\frac{\lambda_{2}(-1)}{\lambda_{1}(-1)}=\eta.$
Then $$f(0,1,k)=\lambda_{2}(1)\left(\eta i^{k}+(-i)^{k}\right),\
f(0,-1,k)=\lambda_{1}(-1)\left(i^{k}+\eta(-i)^{k}\right)
\eqno(3.13)$$ In this condition, we can get that f(2,0,0)$\not=0$.
In fact, since
$[[t_{1}^{1},t_{2}^{-1}],t_{1}^{1}].v_{k}=\frac{q^{-1}-1}{1+q^{-1}}[t_{1}^{2},t_{2}^{-1}].v_{k},$
 we have
 $$\lambda^{2}\left((f(0,-1,k+1)-f(0,-1,k+2))-(f(0,-1,k)-f(0,-1,k+1))\right)$$
$$=\frac{q^{-1}-1}{1+q^{-1}}f(2,0,0)\left(f(0,-1,k)-f(0,-1,k+2)\right),
\eqno(3.14)$$ by (3.4). Substituting (3.13) into (3.14), we have
that
$$\lambda^{2}f(0,-1,k+1)=\frac{q^{-1}-1}{1+q^{-1}}f(2,0,0)f(0,-1,k),\;\forall\ k\in{\bf Z},$$
which implies f(2,0,0)$\not=0$. Since
$[[t_{1}^{2},t_{2}^{1}],t_{2}^{-1}].v_{k}=(1-q^{2})(1-q^{-2})t_{1}^{2}.v_{k}$,
 we can obtain
$$(1-q^{2})(1-q^{-2})=\left(f(0,1,k)-f(0,1,k+2)\right)\left(f(0,-1,k)-f(0,-1,k+2)\right), \eqno(3.15)$$
by (3.4). Substituting (3.11) into (3.15), we have that
$$(1-q^{2})(1-q^{-2})=4(-1)^{k}(\lambda_{1}(1)\lambda_{1}(-1)+\lambda_{2}(1)\lambda_{2}(-1))
+4(\lambda_{1}(1)\lambda_{2}(-1)+\lambda_{2}(1)\lambda_{1}(-1)),$$
for all $k\in {\bf Z}$. Thus
$\lambda_{1}(1)\lambda_{1}(-1)+\lambda_{2}(1)\lambda_{2}(-1)=0$.
This together with
$\lambda_{1}(1)\lambda_{1}(-1)=\lambda_{2}(1)\lambda_{2}(-1)$, we
have
$\lambda_{1}(1)\lambda_{1}(-1)=\lambda_{2}(1)\lambda_{2}(-1)=0$,
which is a contradiction. Therefore, we have proved that if $|x|=1$
then
$\lambda_{1}(1)\lambda_{1}(-1)=\lambda_{2}(1)\lambda_{2}(-1)=0$.

\par Subcase 2. $|x|\neq1$.
\par Without loss of generality, we can suppose $|x|>1$. Then $\lim\limits_{k\rightarrow
+\infty}|x|^{k}=\infty,\ \lim\limits_{k\rightarrow
-\infty}|x|^{k}=0,$ which implies $\lim\limits_{k\rightarrow
+\infty}x^{k}=\infty,\ \lim\limits_{k\rightarrow -\infty}x^{k}=0.$
If $\lambda_{1}(-1)\lambda_{1}(1)\neq0$, we get that
$$\lim\limits_{k\rightarrow
+\infty}\Big(\lambda_{1}(-1)\lambda_{1}(1)(1-2x+x^{2})x^{2k}+
 \lambda_{2}(-1)\lambda_{2}(1)(1-2x^{-1}+x^{-2})x^{-2k}$$
$$\hspace*{\fill} -\left(\lambda_{1}(-1)\lambda_{2}(1)
+\lambda_{2}(-1)\lambda_{1}(1)\right)(x+x^{-1}-2)\Big)=\infty,$$ but
$\lim\limits_{k\rightarrow +\infty}(1-q)(1-q^{-1})$ is a scalar,
which is a contradiction to (3.10). Hence,
$\lambda_{1}(-1)\lambda_{1}(1)=0$.
\par One can deduce that $\lambda_{2}(-1)\lambda_{2}(1)=0$ similarly.
\par All in all, the result that
$\lambda_{1}(-1)\lambda_{1}(1)=\lambda_{2}(-1)\lambda_{2}(1)=0$ is
obtained.
\par Applying the result above to (3.10), we obtain
$$(1-q)(1-q^{-1})=\left(\lambda_{1}(-1)\lambda_{2}(1)
+\lambda_{2}(-1)\lambda_{1}(1)\right)(2-x-x^{-1}).\eqno(3.16)$$ If
$\lambda_{1}(1)=\lambda_{2}(1)=0$ or
$\lambda_{1}(-1)=\lambda_{2}(-1)=0$, substituting this result to
(3.16), we have $(1-q)(1-q^{-1})=0$, which is absurd. Thus $$\bigg
\{
\begin{array}{l}
\lambda_{1}(1)=\lambda_{2}(-1)=0,\\
\lambda_{1}(-1)\neq0,\lambda_{2}(1)\neq 0
\end{array}
  \ \mbox{or}\ \bigg\{
\begin{array}{l}
\lambda_{1}(-1)=\lambda_{2}(1)=0,\\
\lambda_{1}(1)\neq 0,\lambda_{2}(-1)\neq 0.
\end{array}
$$
Without loss of generality, we can assume that
$\lambda_{1}(-1)=\lambda_{2}(1)=0,\ \lambda_{1}(1)\neq 0\ and \
\lambda_{2}(-1)\neq 0.$ Substituting it into (3.16) and (3.9)
respectively, we have
$$(1-q)(1-q^{-1})=\lambda_{2}(-1)\lambda_{1}(1)\left(2-x-x^{-1}\right),\eqno (3.17)$$
$$f(0,1,k)=\lambda_{1}(1)x^{k};\ f(0,-1,k)=\lambda_{2}(-1)x^{-k}. \eqno (3.18)$$
 Setting $\lambda_{1}(1)=a$, from (3.17)
and (3.18) we obtain
$$f(0,1,k)=ax^{k};\ f(0,-1,k)=\frac{(1-q)(1-q^{-1})}{a(2-x-x^{-1})}x^{-k},\
 \forall \ k\in {\bf Z}.$$
\par In a word, we have proved the claim that $$f(0,1,k)=ab^{k}\ \mbox{and}\
f(0,-1,k)=\frac{(1-q)(1-q^{-1})}{a(2-b-b^{-1})}b^{-k}\eqno (3.19)$$
for some $a\in C^{*}\ \mbox{where}\ b^{\pm1}\in{\bf
C}^*\setminus\{1\}$ are roots of the equation (3.6).
\par Considering that $f(\pm1,0,k)=\lambda$ for any
$k\in {\bf Z}$, by using (3.4) and the equation
$$[[t_{1}^{m},t_{2}^{1}],t_{1}^{1}].v_{0}=\frac{(1-q^{m})(q-1)}{1-q^{m+1}}[t_{1}^{m+1},t_{2}^{1}].
v_{0},$$
where $ m\neq 0,\ -1$, we have
$$\frac{(1-q^{m})(q-1)}{1-q^{m+1}}\left(f(0,1,0)-f(0,1,m+1)\right)f(m+1,0,0)$$
$$=\lambda\left(f(0,1,1)-f(0,1,m+1)-f(0,1,0)+f(0,1,m)\right)f(m,0,0).\eqno(3.20)$$
Substituting (3.19) into (3.20), we have that
$$\frac{1-b^{m+1}}{1-q^{m+1}}f(m+1,0,0)=\frac{\lambda(1-b)}{1-q}\frac{1-b^{m}}{1-q^{m}}f(m,0,0)
\;\;\mbox{for\ all}\ m\neq 0,-1. \eqno(3.21)$$ Thus $b^m\neq 1$ for
all $m\in{\bf Z}^*$ since $f(\pm 1,0,0)=\lambda\neq 0$ and $b\neq
1$. ({\bf In fact, one can easily see that Case I discussed in page
7 would not occur from this result.}) Therefore, by the equation
(3.21), we deduce that
$$f(m,0,0)=\left\{\begin{array}{ll}
\left(\frac{\lambda(1-b)}{1-q}\right)^{m}\frac{1-q^{m}}{1-b^{m}},\ \ \ \ \ \ \ \ m\geq1;\\
\frac{q}{b}\left(\frac{\lambda(1-b)}{1-q}\right)^{m+2}\frac{1-q^{m}}{1-b^{m}},\
\ \ m\leq-1.
\end{array} \right. \eqno (3.22)$$
From (3.4) and (3.22) we get
$$f(m,0,k)=\left\{\begin{array}{ll}
\left(\frac{\lambda(1-b)}{1-q}\right)^{m}\frac{1-q^{m}}{1-b^{m}},\ \ \ \ \ \ \ \ m\geq1;\\
\frac{q}{b}\left(\frac{\lambda(1-b)}{1-q}\right)^{m+2}\frac{1-q^{m}}{1-b^{m}},\
\ \ m\leq-1,
\end{array} \right. \eqno (3.23)$$
for any $k\in {\bf Z}$.
\par Applying (3.19) to the equation
$$(1-q^{m})(t_{1}^{m}t_{2}^{j}).v_{k}
=[t_{1}^{m}t_{2}^{j-1},t_{2}^{1}].v_{k},$$
where $m\neq 0,\ j\neq
1$, we have
$$f(m,j,k)=\frac{ab^{k}(1-b^{m})}{1-q^{m}}f(m,j-1,k),\eqno(3.24)$$
where $m\neq0,\ j\neq 1$. And the equation (3.24) implies
$$f(m,j,k)=\left\{\begin{array}{ll}
\left(\frac{ab^{k}(1-b^{m})}{1-q^{m}}\right)^{j-1}f(m,1,k),\ \ \
\ \ \ \ \ \ \ j\geq1,\\
\left(\frac{ab^{k}(1-b^{m})}{1-q^{m}}\right)^{j}f(m,0,k),\ \ \ \ \
\ \ \ \ \ j\leq0,
\end{array} \right.\eqno(3.25)$$
where $m\neq0$. Using (3.19) and $[t_{1}^{m},t_{2}^{1}].
v_{k}=(1-q^{m})(t_{1}^{m}t_{2}^{1}). v_{k}$, we obtain that
$$f(m,1,k)=\frac{ab^{k}(1-b^{m})}{1-q^{m}}f(m,0,k)\ \mbox{for\ all}\ m\in{\bf Z}^*.\eqno(3.26)$$
Substituting (3.26) into (3.25) and considering f(m,0,k)=f(m,0,0)
for all $k\in {\bf Z}$, we have
$$f(m,j,k)=\left(\frac{ab^{k}(1-b^{m})}{1-q^{m}}\right)^{j}f(m,0,0). \eqno (3.27)$$
\par Since $[t_{1}^{1}t_{2}^{j},t_{1}^{-1}]. v_{k}=(q^{-j}-1)t_{2}^{j}.
v_{k}$, we deduce that
$$\lambda\left(f(1,j,k-1)-f(1,j,k)\right)=(q^{-j}-1)f(0,j,k). \eqno(3.28)$$
Applying (3.27) to (3.28), we have
$$f(0,j,k)=\frac{\lambda^{2}(1-b^{j})}{q^{-j}-1}
\left(\frac{ab^{k-1}(1-b)}{1-q}\right)^{j},\ \forall\ j\in{\bf
Z}^*.\eqno(3.29)$$
Using (3.27) and (3.29), we obtain
$$f(m,j,k)=\left\{\begin{array}{ll}
\frac{\lambda^{2}(1-b^{j})}{q^{-j}-1}\left(\frac{ab^{k-1}(1-b)}{1-q}\right)^{j},\
\ \ \ \ \mbox{where}\ m=0,\  j\in{\bf Z}^*;\\
\left(\frac{ab^{k}(1-b^{m})}{1-q^{m}}\right)^{j}f(m,0,0),\ \ \ \ \ \
\ \forall\ m\in{\bf Z}^{*},\ j,k\in {\bf Z},
\end{array} \right. \eqno(3.30)$$
where f(m,0,0) is given by (3.22).
\par From the equation $[t_{1}^{h}t_{2}^{j},t_{1}^{m}t_{2}^{n}].
v_{k}=(q^{jm}-q^{nh})(t_{1}^{h+m}t_{2}^{j+n}). v_{k}$ we can obtain
$$f(m,n,k)f(h,j,k+m)-f(h,j,k)f(m,n,h+k)=(q^{jm}-q^{nh})f(h+m,j+n,k).\eqno
(3.31)$$ Setting m=h=1, j=0 and n=3 in (3.31) and using (3.22) and
(3.30), we obtain
$$\frac{(1+b)^{2}}{b}=\frac{(1+q)^{2}}{q}, \eqno (3.32)$$ which implies b=q or
$q^{-1}$. Thus the equation (3.6) has different roots q and $q^{-1}$
by (3.19), which implies
$q+q^{-1}=\frac{2\lambda^{2}+(1-q)(q^{-1}-1)}{\lambda^{2}}$.
Therefore, we have $\lambda=\pm1$.
\par Substituting b=q and $\lambda=1$ into (3.22) and (3.30), we obtain
$$f(m,j,k)=(aq^{k})^{j},\ \ \ \forall\
(m,j)\in {\bf Z}^{2}\backslash \{(0,0)\}, \ k\in {\bf
Z},\eqno(3.33)$$ which deduces that the $L$-module $V$ be isomorphic
to $V(a,I)$.
\par Substituting b=q and $\lambda=-1$ into
(3.22) and (3.30), we obtain
$$f(m,j,k)=(-1)^{m}(aq^{k})^{j},\ \ \forall\
(m,j)\in {\bf Z}^{2}\backslash \{(0,0)\}, \ k\in {\bf
Z}.\eqno(3.34)$$ Thus $V$ be isomorphic to $V(a,II)$.
\par Substituting $b=q^{-1}$ and $\lambda=1$ into (3.22) and (3.30),
we obtain
$$f(m,j,k)=(-1)^{m+j+1}(aq^{-k-m})^{j},\ \  \forall\
(m,j)\in {\bf Z}^{2}\backslash \{(0,0)\}, \ k\in {\bf
Z},\eqno(3.35)$$ which deduces that $V$ be isomorphic to $V(a,III)$.
\par Substituting $b=q^{-1}$ and $\lambda=-1$ into (3.22) and (3.30),
we obtain
$$f(m,j,k)=(-1)^{j+1}(aq^{-k-m})^{j},\ \ \forall\
(m,j)\in {\bf Z}^{2}\backslash \{(0,0)\}, \ k\in {\bf
Z}.\eqno(3.36)$$ Therefore, $V$ be isomorphic to $V(a,IV)$.
\par In a word, $L$-module $V$ is isomorphic to one of the four classes of {\bf
Z}-graded $L$-modules constructed in Proposition
2.1.$\hspace*{\fill} \blacksquare$

\par From Lemma 3.1 and Lemma 3.2, we obtain our main
result in this paper.
\par{\bf  Theorem 3.3 }$ :$ If $V$ is a {\bf Z}-graded
$L$-module of the intermediate series, then $V$ is isomorphic to
$V(a,I),\ V(a,II),\ V(a,III),\ V(a,IV)$, or the direct sum of some
trivial $L$-modules. $\hspace*{\fill}\blacksquare$

\end{document}